\documentclass[preprint,3p,times]{elsarticle}

\usepackage{amssymb}
\usepackage{indentfirst}
\usepackage{array} 		
\usepackage{newlfont}
\usepackage{microtype}
\usepackage{amscd}

\usepackage[fleqn,tbtags]{amsmath}
\usepackage{amssymb,mathrsfs} 
\usepackage{dsfont}     
\usepackage{mathtools}
\numberwithin{equation}{section}      

\usepackage{latexsym}
\usepackage{amsthm}
\usepackage{fancyhdr}
\usepackage[english,francais]{babel}
\usepackage[applemac]{inputenc}    
\usepackage[T1]{fontenc}

\theoremstyle{definition}                    
\newtheorem{thm}{Theorem}[section]      
\newtheorem{prop}[thm]{Proposition}    
\theoremstyle{definition}               
\newtheorem{dfn}[thm]{Definition}
\newtheorem{ese}[thm]{Example}        
\newtheorem{rem}[thm]{Remark}          

\def\R{\mathbb R}

\def\P{\mathbb P}
\def\1{\mathds{1}}

\biboptions{longnamesfirst,comma}

\journal{}

\begin{document}
\begin{frontmatter}



\title{Clark-Ocone type formula for non-semimartingales with finite quadratic variation\\
Formule de Clark-Ocone generalis\'ee pour non-semimartingales \`a 
variation quadratique finie}
\author[luiss,ensta]{Cristina DI GIROLAMI}
\ead{cdigirolami@luiss.it}
\author[ensta,inria]{Francesco RUSSO}
\ead{francesco.russo@ensta-paristech.fr}
\address[luiss]{Luiss Guido Carli - Libera
    Universit\`a Internazionale degli Studi Sociali Guido Carli di Roma.}
\address[ensta]{ENSTA ParisTech,
Unit\'e de Math\'ematiques appliqu\'ees,
32, Boulevard Victor,
F-75739 Paris Cedex 15 (France)}
\address[inria]{INRIA Rocquencourt 
and Cermics Ecole des Ponts, Projet MATHFI. Domaine de Voluceau, BP 105
F-78153 Le Chesnay Cedex (France).}

\selectlanguage{english}
\begin{abstract}
We provide a suitable framework for the concept of 
finite quadratic variation for processes with values in  a separable Banach space $B$  
using the language  of stochastic calculus via regularizations,
introduced in the  case $B= \R$ by the  second author and P. Vallois. 
To a real continuous process 
$X$ we associate the Banach valued process $X(\cdot)$, called 
{\it window} process, which describes the evolution of $X$ taking into account a memory $\tau>0$.
The natural state space for $X(\cdot)$ is the Banach space of continuous functions on $[-\tau,0]$.
If $X$ is a real finite quadratic variation process, an appropriated It\^o formula is  presented,
from which we derive a generalized Clark-Ocone formula for
non-semimartingales having the
same quadratic variation as Brownian motion.
The representation is based  on solutions of
an infinite dimensional PDE.\\

\selectlanguage{francais}
\noindent
\textbf{R\'esum\'e}\\
\noindent
Nous pr\'esentons un cadre ad\'equat pour le concept de 
variation quadratique finie lorsque le processus de r\'ef\'erence est \`a valeurs
dans un espace  de Banach s\'eparable $B$.
Le langage utilis\'e est celui de l'int\'egrale
via r\'egularisations introduit dans le cas r\'eel
par le second auteur et P. Vallois. 
\`A un processus r\'eel continu $X$, nous associons le processus $X(\cdot)$, appel\'e processus {\it fen\^etre}, qui \`a l'instant $t$, garde 
en m\'emoire le pass\'e jusqu'\`a $t-\tau$. L'espace naturel d'\'evolution pour $X(\cdot)$ est 
l'espace de Banach $B$ des fonctions continues d\'efinies sur $[-\tau,0]$.
Si $X$ est un processus r\'eel \`a variation quadratique finie, 
nous \'enon\c{c}ons une formule d'It\^o appropri\'ee de laquelle nous d\'eduisons une formule de 
Clark-Ocone relative \`a des non-semimartingales
r\'eelles ayant la m\^eme  variation quadratique que le 
mouvement brownien. La repr\'esentation est bas\'ee
sur des solutions d'une EDP infini-dimensionnelle.

\begin{keyword} 
Calculus via regularization\sep Infinite dimensional 
analysis 
\sep Clark-Ocone formula
\sep It\^o formula \sep Quadratic variation \sep Hedging theory without semimartingales.
\MSC[2010]  60H05\sep 60H07\sep  60H30\sep 91G80.
\end{keyword}
\end{abstract}

\end{frontmatter}

\textbf{Version fran\c{c}aise abr\'eg\'ee}\\

\medskip
Dans cette Note nous d\'eveloppons un calcul stochastique
via r\'egularisation de type progressif ({\it forward})
lorsque le processus int\'egrateur $\mathbb{X}$ est \`a valeurs
dans un espace de Banach s\'eparable $B$. Ceci est bas\'e sur
une notion sophistiqu\'ee de {\it variation quadratique}
que nous appellerons  $\chi$-variation quadratique,
o\`u le symbole $\chi$ correspond \`a un sous-espace
$\chi$ du dual du produit tensoriel projectif $B \hat{\otimes}_{\pi} B$.
Le  calcul via  r\'egularisation a \'et\'e introduit
lorsque  $B = \R$ dans \cite{rv} et depuis il a \'et\'e 
\'etudi\'e par de nombreux auteurs qui ont fait avancer  la th\'eorie
et ont produit plusieurs applications.
Le lecteur peut consulter   \cite{Rus05} pour  une  revue
incluant  une liste assez compl\`ete  de r\'ef\'erences.
Dans ce contexte, les auteurs introduisent une
notion de covariation entre  deux processus 
r\'eels $X$ et $Y$, not\'ee $[X,Y]$ qui g\'en\'eralise
le crochet droit usuel lorsque  $X$ et $Y$ sont des semimartingales.
Un vecteur de processus $ \underline X = (X^1, \ldots, X^n)$  
est dit admettre tous ses crochets mutuels si $[X^i,X^j]$
existe pour tous entiers $1 \le i,j \le n$. \\
 Lorsque  $B = \R^{n}$, $\mathbb{X}$ poss\`ede une
  $\chi$-variation quadratique
avec $\chi = (B \hat{\otimes}_{\pi} B)^*$ si et seulement si $\mathbb{X}$ admet tous
ses crochets mutuels. 
On peut voir qu'un processus \`a valeurs dans un  espace  de Banach
  {\it localement semi-sommable}
 $\mathbb{X}$ au sens de  \cite{dincuvisi}, admet une
 $\chi$-variation quadratique avec
 $\chi = (B \hat{\otimes}_{\pi} B)^*$. Dans ce travail nous tra\c{c}ons 
 une \'ebauche  du calcul stochastique via la formule d'It\^{o} 
\'enonc\'ee au Th\'eor\`eme \ref{thm ITONOM}.
Une attention sp\'eciale est consacr\'ee au cas o\`u 
 $B$ est l'espace $C([-\tau,0])$ des fonctions continues d\'efinies sur
 $[-\tau,0]$, pour un certain $\tau >0$, qui est typiquement
un espace de Banach non-r\'eflexif, et \`a une 
formule de Clark-Ocone g\'en\'eralis\'ee.
Soit $T>0$; tout processus r\'eel continu
 $X=(X_{t})_{t\in [0,T]}$ est prolong\'e par continuit\'e
pour $t\notin [0,T]$.\\
Soit $0<\tau\leq T$ et $X$ un processus r\'eel continu; nous  appelons \textbf{fen\^etre} 
le processus \`a  valeurs dans $C([-\tau,0])$ d\'efini par
\begin{equation*}
X(\cdot)=\big(X_{t}(\cdot)\big)_{t\in [0,T]}=\{X_{t}(u):=X_{t+u}; u\in [-\tau,0], t\in [0,T]\} \; .
\end{equation*}
La th\'eorie de l'int\'egration infini-dimensionnelle par
rapport \`a des martingales (ou des semimartingales, 
\cite{dpz, mp, dincuvisi}) n'est pas appliquable, m\^eme
lorsque l'int\'egrateur est la fen\^etre $W(\cdot)$ associ\'ee au mouvement  brownien standard $W$. 
Au-del\`a des difficult\'es qui viennent du fait que 
$C([-\tau, 0]$ n'est pas r\'eflexif, $W(\cdot)$ n'est d'aucune mani\`ere une semimartingale
\`a valeurs dans $C([-\tau,0])$.\\
%
%
Motiv\'es par des applications  li\'ees \`a  la couverture
d'options d\'ependant  de  toute la  trajectoire, nous discutons
une  formule  de type Clark-Ocone visant  \`a  d\'ecomposer 
une classe significative  $h$ de v.a. 
d\'ependant de la trajectoire  d'un processus
 $X$ dont la  variation quadratique vaut
$[X]_t = t$.
Cette formule g\'en\'eralise des r\'esultats inclus dans 
 \cite{schoklo,VaSoBen,crnsm} visant  \`a d\'eterminer des
formules de valorisation et de couverture d'options
vanille o\`u asiatique dans un mod\`ele de prix d'actif
ayant la m\^eme variation quadratique que le mod\`ele de
Black-Scholes.
Si le bruit dans un environnement  stochastique
est  mod\'elis\'e par la d\'eriv\'ee d'un mouvement brownien $W$, le th\'eor\`eme de 
repr\'esentation des martingales
et la formule classique  de Clark-Ocone sont deux outils fondamentaux de calcul. 
Le th\'eor\`eme \ref{cor GHY} et les consid\'erations \`a la fin de la section \ref{sec:CO} montrent que dans 
une certaine mesure une formule de type Clark-Ocone reste valable lorsque la
loi du processus soujacent n'est plus la  mesure de Wiener
mais le processus conserve la m\^eme variation
quadratique que $W$.
Il est en fait possible de repr\'esenter
des  variables al\'eatoires 
$h=H(X_{T}(\cdot))$, o\`u $H:C([-T,0])\longrightarrow \R$,
comme  
\begin{equation}  \label{eq BNV}
h=H_{0}+\int_{0}^{T}\xi_{t}d^{-}X_{t}
\end{equation} 
sous des conditions suffisantes raisonnables sur la 
fonctionnelle  $H$, o\`u
 $H_{0}$ est un nombre r\'eel et 
 $\xi$ est un processus adapt\'e \`a la filtration associ\'ee \`a $X$
qui sont  donn\'es de fa\c{c}on quasi-explicite. Ici $d^{-}X_{s}$
symbolise l' int\'egration progressive (``forward'') via r\'egularisations d\'efinie dans \cite{Rus05}.
Ces quantit\'es sont exprim\'ees \`a l'aide d'une fonctionnelle 
$u:[0,T]\times C([-T,0])\longrightarrow \R$ de classe  
$C^{1,2}\left([0,T[\times C([-T,0]) \right)$ 
qui est solution d'une \'equation aux d\'eriv\'ees partielles; 
la r\'epresentation \eqref{eq BNV} de $h$ a lieu avec $H_{0}=u(0,X_{0}(\cdot))$ et 
 $\xi_{t}=D^{\delta_{0}}u\, (t,X_{t}(\cdot))$, o\`u
 $D^{\delta_{0}}u\,(t,\eta):=D u\,(t,\eta)(\{0\})$, $Du$ symbolisant la 
 d\'eriv\'ee de Fr\'echet par rapport \`a $\eta\in C([-T,0])$; 
 $Du\,(t, \eta)$ est donc une mesure sign\'ee finie.\\
Si $X$ est un mouvement brownien standard $W$ et $h\in \mathbb{D}^{1,2}$, 
l'expression \eqref{eq BNV} coincide avec la formule de Clark-Ocone classique. 

\selectlanguage{english}
\section{Introduction}

In the whole paper $(\Omega,\mathcal{F},\P)$ is a fixed probability space, equipped with
a given filtration $\mathbb{F}=(\mathcal{F}_{t})_{t\geq 0}$ fulfilling the usual conditions, 
$B$ will be a separable Banach space and $\mathbb{X}$ a $B$-valued process. 
If $K$  is a compact set, ${\cal M}(K)$ will denote the space of  Borel 
(signed) measures on $K$.
$C([-\tau,0])$ will denote the space of continuous functions defined
on $[-\tau,0]$ whose topological dual space is ${\cal M}([-\tau,0]) $.
$W$ will always denote an $(\mathcal{F}_{t})$-real Brownian motion. Let $T>0$ be a fixed maturity time. 
All the processes $X=(X_{t})_{t\in [0,T]}$ 
are prolongated by continuity for $t\notin [0,T]$ setting $X_{t}=X_{0}$ for $t\leq 0$ and $X_{t}=X_{T}$ for $t\geq T$. 

We first recall the basic concepts of forward integral and covariation and some one-dimensional 
results concerning calculus via regularization, a fairly complete survey on the subject being \cite{Rus05}. 
For simplicity, all the considered integrator processes will be continuous. 
\begin{dfn}
Let $X$ (respectively $Y$) be a continuous (resp. locally integrable) process.\\
The \textbf{forward integral of $Y$ with respect to $X$} (resp. the \textbf{covariation of $X$ and $Y$}), whenever it exists, is defined as
\begin{equation}\label{def fwd int}
\int_{0}^{t}Y_{s}d^{-}X_{s}:=\lim_{\epsilon\rightarrow 0^{+}}\int_{0}^{t} Y_{s}\frac{X_{s+\epsilon}-X_{s}}{\epsilon}ds  \quad \textrm{in probability for all $t\in [0,T]$} \; ,
\end{equation}
\begin{equation}				\label{eq def cov}				
\left( \textrm{resp.} \left[X,Y\right]_{t}  =\lim_{\epsilon\rightarrow 0^{+}} \frac{1}{\epsilon} \int_{0}^{t} (X_{s+\epsilon}-X_{s})(Y_{s+\epsilon}-Y_{s})ds   \quad \textrm{in the ucp sense with respect to $t$}  \right),\; 
\end{equation}
provided that the limiting process admits a continuous version. 
If $\int_{0}^{t}Y_{s}d^{-}X_{s}$ exists for any $0\leq t< T$; $\int_{0}^{T}Y_{s}d^{-}X_{s}$ will symbolize the 
\textbf{improper forward integral} defined by $\lim_{t\rightarrow T}\int_{0}^{t}Y_{s}d^{-}X_{s} $, whenever it exists in probability.\\
If $[X,X]$ exists then $X$ is said to be a {\bf finite quadratic variation} process. $[X,X]$ will also be denoted by $[X]$ and it will be called {\bf quadratic variation of $X$}. 
If $[X]=0$, then $X$ is said to be a {\bf zero quadratic variation process}. 
If $\mathbb{X}=(X^{1},\ldots, X^{n})$ is a vector of continuous processes we say that it has all its \textbf{mutual covariations} (brackets) if $[X^{i},X^{j}]$ exists for any 
$1\leq i,j\leq n$.
\end{dfn}  
When $X$ is a (continuous) semimartingale (resp. Brownian motion) and $Y$ is an adapted
cadlag  process (resp. such that $\int_0^T Y_s^2 ds <  \infty$ a.s.), the integral $\int_0^\cdot Y_{s} d^-X_{s}$ 
exists and coincides with classical It\^o's integral
$\int_0^\cdot Y_{s} dX_{s}$, see Proposition 6 in \cite{Rus05}. Stochastic calculus via regularization is a  theory
which allows, in many specific cases to manipulate those integrals
when $Y$  is anticipating or $X$ is not a semimartingale.
If $X$, $Y$ are $(\mathcal{F}_{t})$-semimartingales then $[X,Y]$ coincides with the 
classical bracket $\langle X,Y\rangle$, see Corollary 2 in \cite{Rus05}. 
Finite quadratic variation processes will play a central role in this note: this class includes of course all $(\mathcal{F}_{t})$-semimartingales. 
However that class is much richer.
Typical examples of finite quadratic variation processes
are $(\mathcal{F}_t)$-Dirichlet processes. $D$ is called 
$(\mathcal{F}_{t})$-Dirichlet   process if it admits a decomposition $D=M+A$ where $M$ is an
$(\mathcal{F}_{t})$-local martingale and $A$ is a zero quadratic
variation process. It holds in that case $[D]=[M]$. This class of processes generalizes the semimartingales since 
a locally bounded variation process has zero quadratic variation.
A  slight generalization of that notion is
the notion of weak Dirichlet, which was
introduced in \cite{er2}. $X$ is called $(\mathcal{F}_{t})$-weak Dirichlet if it
admits a decomposition $X=M+A$ where $M$ is an $(\mathcal{F}_{t})$
local martingale and $A$ is a process such that $[A,N]=0$ for any
continuous $(\mathcal{F}_{t})$ local martingale $N$. An $(\mathcal{F}_{t})$-weak Dirichlet process is not necessarily a finite quadratic variation process.
On the other hand if $A$ has finite quadratic variation then it holds $[X]=[M]+[A]$. 
Another interesting example is the 
bifractional Brownian motion 
$B^{H,K}$ with parameters $H\in (0,1)$ and $K\in (0,1]$ which has finite quadratic variation if and only if $HK\geq 1/2$, see \cite{rtudor}. 
Notice that if $K=1$, then $B^{H,1}$ coincides with a fractional Brownian motion with Hurst parameter $H\in (0,1)$. 
If $HK=1/2$ it holds $[B^{H,K}]=2^{1-K}t$; if $K\neq 1$ this process is not even Dirichlet with respect to its own filtration. 
Other significant examples are the so-called weak $k$-order Brownian motions, for fixed $k\geq 1$, constructed by \cite{follWuYor}, which are in general not Gaussian.
$X$ is a weak $k$-order Brownian motion if for every $0\leq t_{1}\leq \cdots\leq t_{k}< +\infty$, 
$(X_{t_{1}},\cdots,X_{t_{k}})$ is distributed as $(W_{t_{1}},\cdots, W_{t_{k}})$. 
If $k\geq 4$ then $[X]_{t}=t$.\\
One central object of this work will be the generalization to infinite dimensional valued processes 
of the stochastic integral via regularization, see Definition \ref{def integ fwd}. A stochastic calculus for Banach valued martingales was considered 
by \cite{brez, vanNervUMD} and references therein, generalizing the classical stochastic calculus of \cite{dpz, mp, dincuvisi}.\\
We introduce now a particular Banach valued process.
Given $0<\tau\leq T$ and a real continuous process $X$, we will call \textbf{window process} associated with $X$, the 
$C([-\tau,0])$-valued process denoted by $X(\cdot)$ defined as
\begin{equation*}
X(\cdot)=\big(X_{t}(\cdot)\big)_{t\in [0,T]}=\{X_{t}(u):=X_{t+u}; u\in [-\tau,0], t\in [0,T]\} \; .
\end{equation*}
The window process $W(\cdot)$ associated with the classical Brownian motion $W$ will be called \textbf{window Brownian motion}.
We observe that $W(\cdot)$ is not a $B= C([-\tau,0])$-valued semimartingale
even in the (weak) sense that $ \prescript{}{B^{\ast}}{\langle} \mu, W_t(\cdot) \rangle_B$ is a real semimartingale 
for any $\mu \in B^\ast$. 
In fact setting $\mu=\delta_{0}+\delta_{-\tau/2}$, we get 
$$ Y_t:= \prescript{}{\mathcal{M}([-\tau,0])}{\langle} \mu, W_{t}(\cdot)\rangle_{C([-\tau,0])} =
\int_{-\tau}^0 W_t(u) d\mu(u)  =    W_{t} + W_{t- \frac{\tau}{2}}$$
which is not a semimartingale. In fact its canonical filtration is the filtration
$({\cal F}_t)$ associated with $W$. Taking into account Corollary 3.14 of \cite{crwdp}
$Y$ is an $({\cal F}_t)$-weak Dirichlet process with martingale part $W$.
By uniqueness of the decomposition of a weak Dirichlet process (see Proposition
16 of \cite{Rus05}) $Y$ cannot be an $({\cal F}_t)$-semimartingale. \\
Motivated by the necessity of an It\^o formula available also for $B=C([-\tau,0])$-valued processes, we introduce a 
quadratic variation concept which depends on a subspace $\chi$ of the dual of the tensor square of $B$, 
equipped with the projective topology, denoted by 
$(B\hat{\otimes}_{\pi}B)^{\ast}$, see Definition \ref{DChiQV}. We recall the fundamental identification $(B\hat{\otimes}_{\pi} B)^{\ast}\cong \mathcal{B}(B\times B)$, which denotes 
the space of $\R$-valued bounded bilinear forms on $B\times B$. An It\^o formula for processes admitting a $\chi$-quadratic variation is given in Theorem \ref{thm ITONOM}. 
After formulating a theory for $B$-valued processes with general $B$, in Sections \ref{sec:C} and \ref{sec:CO} we fix the attention on window processes setting $B=C([-\tau,0])$. 
Section \ref{sec:C}, in particular Proposition \ref{pr QV123}, is devoted to the evaluation of $\chi$-quadratic variation for windows
 associated with real finite quadratic variation processes. 
 Suppose that $X$ is a real process such that $[X]_{t}=t$.
In Section \ref{sec:CO} we give a representation result for 
a random variable $h:=H(X_{T}(\cdot))$ where $H:C([-T,0])\longrightarrow \R$ is continuous.
That is of the type $h=H_{0}+\int_{0}^{T}\xi_{s}d^{-}X_{s}$, $H_{0}\in \R$ and $\xi$ adapted process where the integral is considered as the forward 
integral defined in \eqref{def fwd int}. 
More precisely $h$ will appear as $u(T, X_{T}(\cdot))$ where 
$u\in C^{1,2}\left( [0,T[\times C([-T,0]);\R \right)\cap C^{0}\left( [0,T]\times C([-T,0]);\R \right)$ 
solves an infinite dimensional 
partial differential equation of type \eqref{eq SYST}. Moreover 
we will get $H_{0}=u(0,X_{0}(\cdot))$ and $\xi_{s}=D^{\delta_{0}}u\, (s, X_{s}(\cdot))$ where $D^{\delta_{0}}u\,(t,\eta):=D u\,(t,\eta)(\{0\})$; 
$Du$ denotes the Fr\'echet derivative with respect to $\eta \in C([-T,0]$ so $Du\,(t,\eta)$ is a signed measure.
%
\section{Notations}

Symbol $\mathscr{C}([0,T])$ denotes the linear space of continuous real processes 
equipped with the ucp (uniformly convergence in probability) topology, $B^{\ast}$ will be the topological dual of the Banach space $B$. 
We introduce now some subspaces of measures that we will frequently use.
Symbol $\mathcal{D}_{0} ([-\tau,0])$ ( resp. $\mathcal{D}_{0,0} ([-\tau,0]^{2})$), shortly $\mathcal{D}_{0,0}$ (resp. $\mathcal{D}_{0,0}$), will denote the one dimensional Hilbert space of the
 multiples of Dirac measure concentrated at $0$ (resp. at $(0,0)$), i.e.
 \begin{equation}			\label{eq-def D0}
\mathcal{D}_{0} ([-\tau,0]):= \{ \mu \in \mathcal{M}([-\tau,0]);\; s.t. \mu(dx)=\lambda \,\delta_{0}(dx) \textrm{ with } \lambda \in \mathbb{R} \}
\end{equation}
( resp. 
\begin{equation}			\label{eq-def D00}
\mathcal{D}_{0,0} ([-\tau,0]^{2}):= \{ \mu \in \mathcal{M}([-\tau,0]^{2});\; s.t. \mu(dx,dy)=\lambda \,\delta_{0}(dx)\delta_{0}(dy) \textrm{ with } \lambda \in \mathbb{R} \} \; ).
\end{equation}
Symbol $Diag([-\tau,0]^{2})$, shortly $Diag$, will denote the subset of $\mathcal{M}([-\tau,0]^{2})$ defined as follows:
\begin{equation} 		
Diag([-\tau,0]^{2}):=\left\{\mu\in \mathcal{M}([-\tau,0]^{2})\, s.t.\,
\mu(dx,dy)=g(x)\delta_{y}(dx)dy;\, g \in L^{\infty}([-\tau,0])
\right\}
\; .
\end{equation}
$Diag([-\tau,0]^{2})$, equipped with the norm $\| \mu \|_{Diag([-\tau,0]^{2})}= \| g\|_{\infty}$, is a Banach space.

\section{Forward integrals in Banach spaces}		\label{sec: cvr}
In this section we introduce an infinite dimensional stochastic integral via regularization. 
In this construction there are two main difficulties.
The integrator is generally not a semimartingale or 
the integrand may be anticipative; 
$B$ is a general separable, not necessarily reflexive, Banach space.
\begin{dfn}		\label{def integ fwd}  		 
Let $(\mathbb{X}_{t})_{t\in[0,T]}$ (respectively $(\mathbb{Y}_{t})_{t\in[0,T]}$) be a
$B$-valued (respectively a $B^{\ast}$-valued) stochastic process. 
We suppose $\mathbb{X}$ to be continuous and $\mathbb{Y}$ to be strongly measurable (in the Bochner sense) 
such that $\int_{0}^{T}\|\mathbb{Y}_{s}\|_{B^{\ast}} ds <+\infty$ a.s.\\
For every fixed $t\in [0,T]$ we define the {\bf definite forward integral of $\mathbb{Y}$ with respect
to $\mathbb{X}$} denoted by $\int_{0}^{t}\prescript{}{B^{\ast}}{\langle} \mathbb{Y}_{s}, d^{-}\mathbb{X}_{s}\rangle_{B}$ as the following limit in probability:
\begin{equation}		\label{def INTFWD}
\int_{0}^{t} \prescript{}{B^{\ast}}{\langle} \mathbb{Y}_{s}, d^{-}\mathbb{X}_{s}\rangle_{B}: =\lim_{\epsilon\rightarrow
0}\int_{0}^{t} \prescript{}{B^{\ast}}{\langle} \mathbb{Y}_{s},\frac{\mathbb{X}_{s+\epsilon}-\mathbb{X}_{s}}{\epsilon}\rangle_{B} ds		\; .
\end{equation}
We say that the {\bf forward stochastic integral of $\mathbb{Y}$ with respect
to $\mathbb{X}$} exists if the process
\[
\left(\int_{0}^{t} \prescript{}{B^{\ast}}{\langle} \mathbb{Y}_{s}, d^{-}\mathbb{X}_{s}\rangle_{B}\right)_{t\in[0,T]}
\]
admits a continuous version. In the sequel indices $B$ and $B^{\ast}$ will  often be omitted.
\end{dfn}
\section{Chi-quadratic variation}		\label{sec: chiqv}
\begin{dfn}\label{DefChi}
A closed linear subspace  $\chi$ of 
$(B\hat{\otimes}_{\pi}B)^{\ast}$, endowed with its own norm,  such that
\begin{equation}   \label{relazione norma chi}
\|\cdot\|_{ (B\hat{\otimes}_{\pi}B)^{\ast}}  \leq \textrm{const} \cdot \|\cdot\|_{\chi}
\end{equation}
will be called a {\bf  Chi-subspace (of $(B\hat{\otimes}_{\pi}B)^{\ast}$}). 
\end{dfn}
%
%
%
%
Let $\chi$ be a Chi-subspace of $(B\hat{\otimes}_{\pi}B)^{\ast}$, $\mathbb{X}$ be a $B$-valued stochastic process and $\epsilon > 0$. 
We denote by $[\mathbb{X}]^{\epsilon}$, the application
\begin{equation}		\label{eq Xepsilon}
[\mathbb{X}]^{\epsilon}:\chi\longrightarrow \mathscr{C}([0,T])
\hspace{1cm} \textrm{defined by}\hspace{1cm}
\phi
\mapsto
\left( \int_{0}^{t} \prescript{}{\chi}{\langle} \phi,
\frac{J\left(  \left(\mathbb{X}_{s+\epsilon}-\mathbb{X}_{s}\right)\otimes^{2}   \right)}{\epsilon} 
\rangle_{\chi^{\ast}} \,ds 
\right)_{t\in [0,T]}
\end{equation}
where $ J: B\hat{\otimes}_{\pi}B\longrightarrow (B\hat{\otimes}_{\pi}B)^{\ast\ast}$
denotes the canonical injection between a space and its bidual.
\begin{rem} 
\item
\begin{enumerate}
\item We recall that $\chi\subset(B\hat{\otimes}_{\pi}B)^{\ast}$ implies 
$(B\hat{\otimes}_{\pi}B)^{\ast\ast}\subset \chi^{\ast}$. 
\item As indicated, $\prescript{}{\chi}{\langle}\cdot,\cdot\rangle_{\chi^{\ast}}$ denotes the duality
between the space $\chi$ and its dual $\chi^{\ast}$; in fact, by assumption,
 $\phi$ is an element of $\chi$ and element $J\left( \left(\mathbb{X}_{s+\epsilon}-\mathbb{X}_{s}\right)\otimes^{2}   \right)$ naturally belongs to $(B\hat{\otimes}_{\pi}B)^{\ast\ast}\subset \chi^{\ast}$. 
\item The real function $s\rightarrow \langle \phi , J((\mathbb{X}_{s+\epsilon}-\mathbb{X}_{s})\otimes^{2})\rangle$ is integrable 
since $|\langle \phi , J((\mathbb{X}_{s+\epsilon}-\mathbb{X}_{s})\otimes^{2}) \rangle  | \leq \textrm{const}   \|\phi\|_{\chi}  \|\mathbb{X}_{s+\epsilon}-\mathbb{X}_{s}  \|^{2}_{B}$.
\item
With a slight abuse of notation, in the sequel, the application $J$ will be omitted. The tensor product 
$\left(X_{s+ \epsilon}-X_{s}\right)\otimes^{2}$ has to be considered as the element $J\left(   \left(\mathbb{X}_{s+\epsilon}-\mathbb{X}_{s}\right)\otimes^{2}   \right)$ which belongs to $\chi^{\ast}$.
\end{enumerate}
\end{rem}
We give now the definition of the $\chi$-quadratic variation of a $B$-valued stochastic process $\mathbb{X}$. 
\begin{dfn} 			\label{DChiQV}   
Let $\chi$ be a separable Chi-subspace of 
$(B\hat{\otimes}_{\pi}B)^{\ast}$ and $\mathbb{X}$ a $B$-valued stochastic process.
We say that $\mathbb{X}$ {\bf admits a} $\chi$-{\bf quadratic variation} if the following
assumptions are fulfilled.
\begin{description}
\item [H1] For every sequence $(\epsilon_{n})\downarrow 0$ there is a
subsequence $(\epsilon_{n_{k}})$ such that 
\[
\sup_{k}\int_{0}^{T}\sup_{\|\phi\|_{\chi}\leq 1}\left|\prescript{}{\chi}{\langle} \phi,\frac{(\mathbb{X}_{s+\epsilon_{n_{k}}}-\mathbb{X}_{s})\otimes^{2}}{\epsilon_{n_{k}}}\rangle_{\chi^{\ast}} \right|ds\quad < +\infty  \; ,\; a.s.
\]
\item [H2] It exists an 
application denoted by $[\mathbb{X}]:\chi\longrightarrow
  \mathscr{C}([0,T])$, such that 
\begin{equation}		\label{convUC}
[\mathbb{X}]^{\epsilon}(\phi)\xrightarrow[\epsilon\longrightarrow 0_{+}]{ucp} [\mathbb{X}](\phi)
\end{equation}
for all $\phi\in \mathcal{S}$, where $\mathcal{S}\subset \chi$ such that $\overline{Span({\cal S})}=\chi$.
\end{description}
\end{dfn}

We formulate a technical proposition which is stated in Corollary 4.38 in \cite{DGR}. Its proof is based on 
Banach-Steinhaus and separability arguments.
\begin{prop}			\label{prBG}
Suppose that $\mathbb{X}$ admits a $\chi$-quadratic variation.
\begin{enumerate}
\item 
Relation \eqref{convUC} holds for any $\phi\in \chi$ and $[\mathbb{X}]$ is a linear continuous application. 
In particular $[\mathbb{X}]$ does not depend on $\mathcal{S}$.
\item
There exists a $\chi^{\ast}$-valued measurable process $(\widetilde{[\mathbb{X}]})_{0\leq t\leq T}$, cadlag and with bounded variation on $[0,T]$ such that 
$\widetilde{[\mathbb{X}]}_{t}(\cdot)(\phi)=[\mathbb{X}](\phi)(\cdot,t) $ a.s. for any $t\in [0,T]$ and $\phi\in \chi$.
\end{enumerate}
\end{prop}

The existence of $\widetilde{[\mathbb{X}]}$ guarantees that 
$[\mathbb{X}]$ admits a proper version which allows to consider it as a 
pathwise integral.

\begin{dfn}
When $\mathbb{X}$ admits a $\chi$-quadratic variation, 
the $\chi^{\ast}$-valued measurable
process  $(\widetilde{[\mathbb{X}]})_{0\leq t\leq T}$ appearing in Proposition 
\ref{prBG}, is called  $\chi$-{\bf quadratic variation} of $\mathbb{X}$. 
Sometimes, with a slight abuse of notation,  
even $[\mathbb{X}]$ will be called $\chi$-quadratic variation
and it will be confused with $\widetilde{[\mathbb{X}]}$.
\end{dfn}

%
%
%
%
%
\begin{dfn} 			\label{DQVWS}
We say that a continuous $B$-valued process $\mathbb{X}$ admits 
{\bf global quadratic variation} if it admits a $\chi$-quadratic variation with
 $\chi=(B\hat{\otimes}_{\pi} B)^{\ast}$. In particular $\widetilde{[\mathbb{X}]}$ takes values ``a priori'' in $(B\hat{\otimes}_{\pi} B)^{\ast\ast}$. 
\end{dfn}
The natural generalization of quadratic variation for a $B$-valued \textit{locally semi summable} process is a $(B\hat{\otimes}_{\pi}B)$-valued process, called 
the \textit{tensor quadratic variation}, as it was introduced by 
\cite{dincuvisi} and \cite{mp}. Unfortunately, the tensor quadratic variation does not exist in several contexts. 
For instance, the window Brownian motion $W(\cdot)$,
which is our fundamental example, does not admit it, see Remark \ref{rem WBM}. 
That notion is related to a strong convergence in $B\hat{\otimes}_{\pi}B$ while our concept of 
global quadratic variation is related to a weak star convergence in its bidual. 
The global quadratic variation generalizes the tensor quadratic variation one:
if $\mathbb{X}$ admits a tensor quadratic variation then it admits a global quadratic variation 
and those quadratic variations are equal, see Section 6.3 in \cite{DGR} for       details. 
When $B$ is the finite dimensional space 
$\R^{n}$, 
$\mathbb{X}$ admits a tensor quadratic variation and if and only if $\mathbb{X}$ admits a global quadratic variation. In that case previous properties are also equivalent to 
the existence of all the mutual brackets in the sense of \cite{Rus05}.\\
%
%
%
%
%
%
%
\section{It\^{o}'s formula}                                      	
\label{sec: ito}
The classical It\^o formulae for stochastic integrators $\mathbb{X}$ with values in an infinite dimensional space appear in Section 4.5 of \cite{dpz} for the Hilbert separable case and in Section 3.7 in \cite{mp}, see also \cite{dincuvisi}, as far as the Banach case is concerned; they involve processes admitting a tensor quadratic variation. 
We state now an It\^{o} formula in the general separable Banach space which do not necessarily have 
a tensor quadratic variation but they have rather a $\chi$-quadratic variation, where $\chi$ is some Chi-subspace where the second order Fr\'echet derivative lives. 
This type of formula is well suited for $C([-\tau,0])$-valued integrators as for instance window processes; this will be developed in Sections \ref{sec:C} and \ref{sec:CO}.
In the sequel if $F:[0,T]\times B\longrightarrow \mathbb{R}$ then (if it exists) 
$DF$ (resp. $D^{2}F$) stands for the first (resp. second) order Fr\'echet derivative 
with respect to the $B$ variable.\\
\begin{thm}  			 \label{thm ITONOM}
Let $B$ be a separable Banach space, $\chi$ be a Chi-subspace of $(B\hat{\otimes}_{\pi}B)^{\ast}$ 
and $\mathbb{X}$ a $B$-valued continuous process
admitting a $\chi$-quadratic variation.
Let $F:[0,T]\times B\longrightarrow \mathbb{R}$ of class $C^{1,2}$ Fr\'echet.
such that
\begin{equation} 		
D^{2}F:[0,T]\times B\longrightarrow \chi \subset
(B\hat{\otimes}_{\pi}B)^{\ast} \textrm{ is continuous with respect to } \chi \; .
\end{equation}
Then the forward integral 
\[
\int_{0}^{t}\prescript{}{B^{\ast}} {\langle} DF(s,\mathbb{X}_{s}),d^{-}\mathbb{X}_{s}\rangle_{B} \; , \quad t\in [0,T] \; ,
\] 
exists and the following formula holds
\begin{equation}					\label{eq ITONOM}
F(t,\mathbb{X}_{t})=F(0, \mathbb{X}_{0})+\int_{0}^{t}\partial_{t}F(s,\mathbb{X}_{s})ds+\int_{0}^{t}\prescript{}{B^{\ast}}{\langle} DF(s,\mathbb{X}_{s}),d^{-}\mathbb{X}_{s}\rangle_{B} 
+\frac{1}{2}\int_{0}^{t} \prescript{}{\chi}{\langle} D^{2}F(s,\mathbb{X}_{s}),
d\widetilde{[\mathbb{X}]}_{s}\rangle_{\chi^{\ast}} \ a.s.
\end{equation}
\end{thm}
Its proof is 
given in Section 8 of \cite{DGR}.

\section{Evaluation of $\chi$-quadratic variations for window processes}		\label{sec:C}

From this section we fix $B$ as the Banach space $C([-\tau,0])$. 
In this section we give some examples of Chi-suspaces and then we give some evaluations of $\chi$-quadratic variations for window processes $\mathbb{X}=X(\cdot)$. 
For illustration of possible applications of It\^o formula \eqref{eq ITONOM}, consider the following functions. Let $H: B\rightarrow \mathbb{R}$ and $\eta \in B$ 
defined by
\begin{equation}\label{eq exam} 
a) \quad H(\eta)=f(\eta (0)), \quad f\in C^{2}(\mathbb{R})      \, ; \quad  
b) \quad H(\eta)=\left( \int_{-\tau}^{0} \eta(s)ds \right)^{2}     \, ; \quad 
c) \quad H(\eta)=\int_{-\tau}^{0} \eta^{2}(s)ds   \, .
\end{equation}
Those functions are of class $C^{2}(B)$; computing the second order Fr\'echet derivative $D^{2}H:B \rightarrow (B\hat{\otimes}_{\pi}B)^{\ast}$ 
we obtain the following:
\begin{equation}
a) \quad D_{dx\, dy}^{2}H(\eta)=f''(\eta(0)) \, \delta_{0}(dx)\delta_{0}(dy)    \, ;  \quad 
b) \quad D^{2}H(\eta)= 2  \mathds{1}_{[-\tau,0]^{2}}			   \, ; \quad 
c) \quad D^{2}_{dx\, dy}H(\eta)=2\, \delta_{x}(dy)dx  \, .
\end{equation} 
In all those examples, $D^{2}H(\eta)$ lives in a particular Chi-subspace $\chi$. Respectively we have $D^{2}H:B\rightarrow \chi$ continuously with
\begin{equation}
a) \quad \chi= \mathcal{D}_{0, 0} ([-\tau,0]^{2}) \, ; \quad 
b) \quad \chi=L^{2}([-\tau,0]^{2})  \, ; 
\quad c) \quad \chi=Diag([-\tau,0]^{2}) \, .
\end{equation}
Other examples of Chi-subspaces are $\mathcal{M}([-\tau,0]^{2})$ and its subspace $\chi^{0}([-\tau,0]^{2})$, (shortly $\chi^{0}$), defined by 
\begin{equation}
\chi^{0}([-\tau,0]^{2}):=( \mathcal{D}_{0}([-\tau,0]) \oplus L^{2}([-\tau,0]))\hat{\otimes}^{2}_{h}  \; ,
\end{equation}
where $\hat{\otimes}_{h}$ stands for the Hilbert tensor product.
The latter one will intervene in Theorem \ref{cor GHY} in relation with the generalized Clark-Ocone formula.
We evaluate now some $\chi$-quadratic variations of window processes.
%
%
\begin{prop}		\label{prop ZQVHC}
Let $X$ be a real valued process with H\"older
continuous paths of parameter $\gamma>1/2$. 
Then $X(\cdot)$ admits a zero global quadratic variation.
\end{prop}

\begin{ese}
Examples of real processes with H\"older
continuous paths of parameter $\gamma>1/2$ are fractional Brownian
motion $B^{H}$ with $H>1/2$ or a bifractional Brownian
motion $B^{H,K}$ with $HK>1/2$.
\end{ese}
\begin{rem} 		\label{rem WBM}	
The window Brownian motion $W(\cdot)$ does not admit a global (and therefore not a tensor) quadratic variation because 
Condition \textbf{H1} is not verified. In fact it is possible to show that 
\begin{equation}			\label{eq QS32}
\int_{0}^{T} \frac{1}{\epsilon} \left\| W_{u+\epsilon}(\cdot)-W_{u}(\cdot) \right\|^{2}_{B}du \geq T \, A^{2}(\tilde{\epsilon}) \ln (1/\tilde{\epsilon})  
\quad \textrm{where} \quad \tilde{\epsilon}=\frac{2\epsilon}{T}
\end{equation}
and $(A(\epsilon))$ is a family of non negative r.v. such that $\lim _{\epsilon\rightarrow 0}A(\epsilon)=1$ a.s.\\
%
\end{rem}
%
%
%
%
\begin{prop}   \label{pr QV123}
Let $X$ be a real continuous process with finite quadratic variation $[X]$ and $0<\tau\leq T$. The following properties hold true. 
\begin{enumerate}
\item [1)] $X(\cdot)$ admits zero $L^{2}([-\tau,0]^{2})$-quadratic variation.
\item [2)] $X(\cdot)$ admits a $\mathcal{D}_{0,0}([-\tau,0]^{2})$-quadratic variation given by 
\begin{equation}		
[X(\cdot)] (\mu)=\mu(\{0,0\})[X], \quad \forall \mu \in \mathcal{D}_{0,0}([-\tau,0]^{2}).
\end{equation}
\item [3)] $X(\cdot) $  admits a $\chi^{0}([-\tau,0]^{2})$-quadratic 
variation which equals 
\begin{equation}		
[X(\cdot)](\mu)=\mu(\{0,0\}) [X], \quad  \forall \mu \in \chi^{0}([-\tau,0]^{2}).
\end{equation}
\item [4)] $X(\cdot) $  admits a $Diag$-quadratic 
variation given by 
\begin{equation}			\label{eq QV DIAG tau}
\mu \mapsto  [X(\cdot)]_{t}(\mu)=\int_{0}^{t\wedge \tau } g(-x) [X]_{t-x} dx \hspace{2cm}		t\in[0,T] \; ,
\end{equation}
where $\mu$ is a generic element in $Diag([-\tau,0]^{2})$ of type $\mu(dx,dy)=g(x)\delta_{y}(dx)dy$, 
with associated $g$ in $L^{\infty}([-\tau,0])$.
\end{enumerate}
\end{prop}
\begin{rem} 
We remark that in the treated cases, the quadratic variation $[X]$ of the 
real finite quadratic variation process $X$ insures the existence of (and completely determines) the $\chi$-quadratic variation. 
For example if $X$ is a real finite quadratic variation process such
 that $[X]_{t}=t$, then 
$X(\cdot)$ has the same $\chi$-quadratic variation as the window 
Brownian motion for the $\chi$ mentioned in the above proposition.
\end{rem}
%


%
%

\section{A generalized Clark-Ocone formula}		\label{sec:CO}		

In this section we will consider $\tau=T$ and we recall that $B=C([-T,0])$. 
Let $X$ be a real  stochastic process such that $X_{0}=0$ and $[X]_{t}=t$. 
Let $H:C([-T,0])\longrightarrow \R$ be a Borel functional; we aim at representing the random variable 
\begin{equation}
h=H(X_{T}(\cdot))\; .
\end{equation}
The main task will consist in looking for classes of functionals 
$H$ for which there is $H_{0}\in \R$ and a predictable process $\xi$ with 
respect to the canonical filtration of $X$ such that $h$
admits the  representation
\begin{equation}	\label{eq 11}
h=H_{0}+\int_{0}^{T}\xi_{s}d^{-}X_{s} \; .
\end{equation} 
Moreover we look for an explicit expression for $H_{0}$ and $\xi$.
%
As a consequence of  
It\^{o}'s formula \eqref{eq ITONOM} for path dependent functionals of the process we will observe that, in those cases,
it is  possible to find a 
function $u$ which solves an infinite dimensional PDE and which gives at the same time the representation result \eqref{eq 11}. 
One possible representation is the following.
%
%
\begin{thm}			\label{cor GHY}
Let $H:C([-T,0])\longrightarrow\R$ be a Borel functional.  
Let $u\in C^{1,2}\left([0,T[\times C([-T,0]) \right)\cap C^{0}\left([0,T]\times C([-T,0]) \right)$ such that 
$x\mapsto D_{x}^{ac}u\,(t,\eta)$ has bounded variation, for any $t\in [0,T]$, $\eta\in C([-T,0])$
and  $D^{ac}u\, (t,\eta)$ is the absolute continuous part of measure $Du\,(t,\eta)$.
We suppose moreover that  $(t, \eta) \mapsto  D^{2}u\, (t,\eta)$
  takes values   in   $\chi^{0}([-T,0]^{2})$ and it is continuous.
Suppose that $u$ is a solution of 
\begin{equation}			\label{eq SYSTOK}
\begin{dcases}
\partial_{t}u\,(t,\eta)+
\int_{]-t,0]}D^{ac}u\,(t,\eta)\,d\eta+\frac{1}{2}D^{2}u\,(t,\eta)(\{0,0\}) =0 & \\
u(T,\eta)=H(\eta) & \\
\end{dcases}
\end{equation}
where the integral $\int_{]-t,0]}D^{ac}u\,(t,\eta)\,d\eta$ has to be understood via an integration by parts as follows:
\[
\int_{]-t,0]}D^{ac}u\,(t,\eta)\,d\eta=D^{ac}u\,(0,\eta)\eta(0)-D^{ac}u\,(-t,\eta)\eta(-t)-\int_{]-t,0]}\eta(x)\,D_{dx}^{ac}u\,(t,\eta)  \; .
\]
Then the random variable $h:=H(X_{T}(\cdot))$ 
admits the following representation
\begin{equation}			\label{eq REPRE}
h=
u(0,X_{0}(\cdot))+\int_{0}^{T}D^{\delta_{0}}u(t,X_{t}(\cdot))d^{-}X_{t}			\; .
\end{equation} 
\end{thm}
\qed

Sections 9.8 and 9.9 in \cite{DGR} provide different reasonable conditions on $H:C([-T,0])\longrightarrow \R$ 
such that there is a function $u$ solving PDE \eqref{eq SYSTOK} in general situations, i.e. fulfilling the hypotheses of Theorem \ref{cor GHY}. 
When $H:C([-T,0])\subset L^{2}([-T,0])\longrightarrow \R$ is $C^{3}\left(L^{2}([-T,0]) \right)$ such that $D^{2}H \in L^{2}([-T,0]^{2})$ with some other 
minor technical conditions, Theorem 9.41 in \cite{DGR} furnishes explicit solutions to \eqref{eq SYSTOK}. Another case for which it is possible to do the same 
is given by Proposition 9.53 in \cite{DGR},
where $h$ depends (not necessarily smoothly) on a finite number of Wiener integrals of the type $\int_{0}^{T}\varphi(s)d^{-}X_{s}$ and $\varphi\in C^{2}(\mathbb{R})$.

\begin{rem} {In relation to Theorem \ref{cor GHY} we observe the following.}
\begin{itemize}
\item Only pathwise considerations intervene and there is no need to suppose that the law of $X$ is Wiener measure.
\item 
Since $H(\eta)=u(T,\eta)$, we observe that $H$ is automatically continuous by hypothesis $u\in C^{0}\left([0,T]\times C([-T,0]) \right)$.
\item Let us suppose $X=W$.
\begin{enumerate}
\item Making use of probabilistic technology, \eqref{eq REPRE} holds in some cases even if $H$ 
is not continuous and $h\notin L^{1}(\Omega)$; we refer to Section 9.6 in \cite{DGR} for this type of results.
\item If $\int_{0}^{T}\xi_{s}^{2}ds<+\infty$ a.s., then the forward integral $\int_{0}^{T}\xi_{t}\,d^{-}W_{t}$ coincides with the
It\^o integral $\int_{0}^{T}\xi_{t}\,dW_{t}$.
\item If the r.v. $h=H(W_{T}(\cdot))$ belongs to $\mathbb{D}^{1,2}$, by uniqueness of the martingale representation theorem and point 2., 
we have $H_{0}=\mathbb{E}[h]$ and $\xi_{t}=\mathbb{E}\left[D^{m}_{t} h\vert \mathcal{F}_{t} \right]$, where $D^{m}$ is the Malliavin gradient; this agrees with Clark-Ocone formula.
\end{enumerate}
\item If $X$ is not a Brownian motion, in general $H_{0}\neq \mathbb{E}[h]$ since $\mathbb{E}\left[ \int_{0}^{T}\xi_{t}d^{-}X_{t}\right]$ 
does not generally vanish. In fact $\mathbb{E}[h]$ will specifically depend on the unknown law of $X$.
\end{itemize}
\end{rem}
\begin{rem}
The assumption $[X]_{t}=t$ is not crucial. 
With some more work it is possible to obtain similar representations even if $[X]_{t}=\int_{0}^{t}a^{2}(s,X_{s})ds$ for a large class of 
$a:[0,T]\times \mathbb{R}\longrightarrow\mathbb{R}$. As a limiting case we show this possibility when $[X]=0$ and 
$h=f\left( \int_{0}^{T}\varphi^{1}(s)d^{-}X_{s}, \ldots, \int_{0}^{T} \varphi^{n}(s)d^{-}X_{s}\right)$ with $\varphi_{i}\in C^{2}([0,T])$ and $f\in C^{2}(\mathbb{R}^{n})$.
We define $V_{t}=u\left( t, X_{t}(\cdot)\right)$ where $u:[0,T]\times C([-T,0])\longrightarrow \R$ defined by
$$
u(t, \eta)=f\left(\int_{]-t,0]}
 \varphi^{1}(s+t)d\eta(s), \ldots, \int_{]-t,0]} \varphi^{n}(s+t)d\eta(s) \right)\; .
$$
After an application of the finite dimensional It\^o formula for finite quadratic variation processes, 
see Proposition 2.4 in \cite{rg2}, and some further calculations, we have 
\begin{equation}	\label{eq ReprH}
h=f(0, \ldots, 0)+\int_{0}^{t} \xi_{s}  d^{-}X_{s}
\end{equation}
with $\xi_{t}= \sum_{i=1}^{n}\partial_{i}f \left(\int_{-t}^{0}\varphi^{1}(s+t)d^{-}X_{s},\ldots, \int_{-t}^{0}\varphi^{n}(s+t)d^{-}X_{s})\right) \varphi^{i}(t)$.
On the other hand, we observe that $u$ solves the PDE $\partial_{t}u + \int_{]-t,0]}D^{ac}u\,(t,\eta)\,d\eta  =0$, 
which is of the same type of \eqref{eq SYSTOK}. 
Representation \eqref{eq ReprH} can be also established via Theorem \ref{thm ITONOM}, 
taking into account that
$X(\cdot)$ admits zero $\chi^{0}$-quadratic variation.
\end{rem}

In chapter 9 in \cite{DGR} we enlarge the discussion presented in Theorem \ref{cor GHY}.
We can give examples where $u:[0,T]\times C([-T,0])\longrightarrow \R$ of class 
$C^{1,2}\left( [0,T[\times C([-T,0]);\R \right)\cap C^{0}\left( [0,T]\times C([-T,0]);\R \right)$
with $D^2 u \in \chi_0$  such that \eqref{eq REPRE} holds 
and $u$ solves an infinite dimensional PDE of the type
\begin{equation}		\label{eq SYST}
\begin{dcases}
\partial_{t}u(t,\eta)+\prescript{ `` }{ }{\int_{]-t,0]}} D^{ac}u\,(t,\eta)\; d\eta^{''} +\frac{1}{2}\langle D^{2}u\,(t,\eta)\; ,\; \1_{D}\rangle=0  & \\
u(T,\eta)  =H(\eta)   &
\end{dcases}
\end{equation}
where 
$
\1_{D}(x,y):=\left\{ 
\begin{array}{ll}
1  & \textrm{if } x=y, \; x,y\in [-T,0]\\
0 & \textrm{otherwise}
\end{array}
\right.
$.
The integral ``
$\int_{]-t,0]}D^{ac}u\,(t,\eta)\; d\eta$'' has to be suitably defined and term 
$\langle D^{2}u\,(t,\eta)\; ,\; \1_{D}\rangle$ indicates the evaluation of the second order derivative on the 
diagonal of the square $[-T,0]^{2}$.\\
We observe that solution of \eqref{eq SYSTOK} are also solutions of \eqref{eq SYST} since 
$\langle D^{2}u\,(t,\eta), \1_{D}\rangle= D^{2}u\,(t,\eta)(\{0,0 \})$ because $D^2 u$ takes values in $\chi^0$.

\bigskip
\bibliographystyle{elsarticle-harv}
\bibliography{biblio}







\end{document}